\documentclass[a4paper,11pt]{amsart}
  \usepackage{amsmath}
  \usepackage{url,graphicx}
  \usepackage{amssymb, color, pstricks-add, manfnt}
  \usepackage{amsmath, amsthm,  amsfonts}
  \usepackage{mathrsfs}

  \theoremstyle{plain}
  \newtheorem{Theorem}{Theorem}[section]

  \newtheorem{Corollary}{Corollary}[section]

    \theoremstyle{remark}
  \newtheorem{remark}{Remark}

  \numberwithin{equation}{section}
  \numberwithin{remark}{section}

\renewcommand{\baselinestretch}{1.00}
\begin{document}

\title{On the maximum principle for general linear elliptic equations.}

\author{Neil S. Trudinger}
\address{Mathematical Sciences Institute, The Australian National University, Canberra ACT 0200, Australia,
School of Mathematics and Applied Statistics, University of Wollongong, Wollongong, NSW 2522, Australia}
\email{Neil.Trudinger@anu.edu.au}

\date{\today}



\thanks{Research supported by Australian Research Council Grants (DP170100929, DP230100499)}

\abstract {We consider maximum principles and related estimates for linear second order elliptic partial differential operators in n-dimensional Euclidean space, which improve previous results, with H-J Kuo, through sharp $L^p$ dependence on the drift coefficient $b$. As in our previous work, the ellipticity is determined through the principal coefficient matrix $\mathcal A$ lying in sub-cones of the positive cone, which are dual cones of the G\aa rding $k$-cones, for $k = 1,\cdots, n$. Our main results are maximum principles for bounded domains, which extend those of Aleksandrov in the case $k=n$, together with extensions to unbounded domains, depending on appropriate integral norms of $\mathcal A$, and corresponding local maximum principles. We also consider applications to local estimates in the uniformly elliptic case, including extensions of the Krylov-Safonov H\"older and Harnack estimates.}

\endabstract

\maketitle

\baselineskip=12.8pt
\parskip=3pt
\renewcommand{\baselinestretch}{1.38}

\section{Introduction}\label{Section 1}
\vskip10pt

In this paper we are concerned with maximum principles for general linear second order partial differential elliptic operators of the form,

\begin{equation}\label{1.1}
Lu = \mathcal A. D^2 u + b.Du = a^{ij}D_{ij}u + b_iD_iu,
\end{equation}

\noindent where $\mathcal A = [a^{ij]}$ is a positive measurable mapping from a domain $\Omega$ in Euclidean $n$-space $\mathbb{R}^n$ into $\mathbb{S}^n$, the linear space of $n\times n$ real symmetric matrices, and $b = (b_1,\cdots b_n)$ is a measurable mapping from $\Omega$ into $\mathbb{R}^n$. In particular, we improve the corresponding results of Kuo and Trudinger in \cite{KT 2007}, with respect to the conditions assumed on the drift coefficient $b$, together with their applications to corresponding improvements of the Pucci conjecture bounds in \cite{T 2020}. 

In accordance with \cite{KT 2007}, the crucial algebraic quantities in our maximum principles are the elementary symmetric polynomials $S_k$ in $\mathbb R^n$ given by 
\begin{equation}\label{1.2}
 = \sum_{1\leq i_1<\cdots<i_k\leq n}\lambda_{i_1}\cdots \lambda_{i_k}, \quad k=1,\cdots,n,
\end{equation}
and their associated cones,
\begin{equation}\label{1.3}
\Gamma_k = \{\lambda \in \mathbb{R}^n\ | \ S_j[\lambda]>0, \ \forall j=1,\cdots,k \},
\end{equation}
and closed dual cones,
\begin{equation}\label{1.4}
\bar\Gamma^*_k = \{\lambda \in \mathbb{R}^n\ | \ \lambda.\mu \ge0, \ \forall \mu\in \Gamma_k\}.
\end{equation}
Normalised dual functions $\rho^*_k$ are defined on $\bar\Gamma^*_k$ as follows:
\begin{equation}\label{1.5}
\rho^*_k(\lambda) = \inf\{ \frac{\lambda.\mu}{n} \ | \ \mu\in \Gamma_k, \ S_k(\mu) \ge  \binom{n}{k} \}.
\end{equation}
As with the corresponding normalised elementary symmetric functions $\rho_k$ on $\bar \Gamma_k$, given by 
$\rho_k=\{S_k(\lambda)/ \binom{n}{k}\}^{1/k}$, the dual functions $\rho^*_k$ are nondecreasing along rays and concave. 

Our  dual cones, corresponding to the open cones $\Gamma_k$, are then given by 
\begin{equation} \label{1.6}
\Gamma^*_k := \{\lambda\in \bar\Gamma^*_k | \rho^*_k(\lambda)>0 \}
\end{equation}

\noindent Note that it is readily shown by approximation that $\bar\Gamma^*_k$ is indeed the closure of $\Gamma^*_k$ and that
from the concavity of $\rho_k$, it also follows that $\Gamma^*_k$ is open in $\mathbb R^n$ for $k > 1$; (see Remark 1.1).

Clearly, $\Gamma_k \subset \Gamma_l$ for $k\ge l$ and $\Gamma_1$ is the half-space $\{\sum \lambda_i > 0\}$, while $\Gamma_n$ is the positive cone. Consequently the dual cones $\Gamma^*_k \subset \Gamma^*_l$  for $k\le l$, with $ \Gamma^*_1$ the open ray through $(1, \ldots, 1)$ and $\Gamma^*_n =  \Gamma_n$.
Note also that our notation here differs from that in \cite{KT 2007}, where we have used $\Gamma^*_k$ to denote the closed cone $\bar\Gamma^*_k$.

Next we say that a symmetric $n\times n$ matrix $A \in \Gamma_k, (\Gamma^*_k, \bar \Gamma_k, \bar \Gamma_k)$ if its eigenvalues $\lambda \in \Gamma_k, (\Gamma^*_k, \bar \Gamma_k, \bar \Gamma_k)$ and define $\rho^*_k (A) = \rho^*_k(\lambda)$. It follows that $\rho^*_k (A)\ge \rho^*_l(A)$ for $k\ge l$, $A\in \Gamma^*_k$ while
$\rho^*_n A) = (\text{det} A)^{1/n}$ and  $\rho^*_1(A)=\lambda_{min}(A)$, where  
$\lambda_{min}(A)$ denotes the minimum eigenvalue of $A$. Expressing the eigenvalues $\lambda(\mathcal A) = (\lambda_1, \dots, \lambda_n)$ in non-decreasing order, so that $\lambda_{min} = \lambda_1$, we also have the following estimate, for $\mathcal A \in \Gamma ^*_k$:

\begin{equation} \label{1.7} 
\rho^*_k(\mathcal A) \le \frac{k}{n} {\binom{n}{k}}^{1/k} (\lambda_1 \dots \lambda_k)^{1/k},
\end{equation}

\noindent with equality in the cases $k=1$ and $k=n$.  To prove \eqref{1.7}, we fix a vector $\mu = (\mu_1, \dots, \mu_n)$ in the definition \eqref{1.5} by setting 

$$ \mu_i = {\binom{n}{k}}^{1/k} (\lambda_1 \dots \lambda_k)^{1/k} \lambda^{-1}_i, $$

\noindent for $i= 1, \dots k$, and $\mu_i =0$ for $i= k+1, \dots n$. 

\begin{remark} 
We remark here that the inequality \eqref{1.7} will also be sharp for $1<k<n$, as equality will hold there along the ray passing through the point where $\lambda_i = 1$, for $i\le k$ and $\lambda_i = k$, for $i>k$. This can be proved by showing that the above choice of $\mu$ will then be the unique critical point for the associated minimisation problem, which is a consequence of the concavity of $\rho_k$. Similarly we can use the uniqueness of an optimal $\mu = \mu(\lambda)$ in \eqref{1.5} for $\lambda\in \Gamma^*_k$ to verify the openness of $\Gamma_k^*$. 
\end{remark}

The following maximum principle, for bounded $\Omega$, now extends that in Theorem 1.1 of \cite{KT 2007} to general linear operators $L$, in the cases, $k\ge n/2$. 

\begin{Theorem} Let $L$ be an operator of the form \eqref{1.1} with coefficients  $\mathcal A (\Omega) \subset \Gamma^*_k$, $b/\rho^*_k(\mathcal A)\in 
L^p(\Omega)$ for $k \ge n/2$ and some   $p>n$, ($\ge n$ if $k=n$). Then if  $u \in W^{2,q}_{loc} (\Omega)\cap C^0(\bar\Omega)$, $(Lu)^- / \rho^*_k(\mathcal A) \in L^q(\Omega)$, for $q\ge k$ if $k>n/2$ and $q>n/2$ if $k=n/2$,  we have the estimate,

\begin{equation} \label{1.8}
\sup_\Omega u \le \sup_{\partial\Omega} u + C (\text{diam}\Omega)^{2-n/q}  ||\frac{(Lu)^-}{\rho^*_k(\mathcal A)}||_{L^q(\Omega)},
\end{equation}

\noindent where $C$ is a constant depending on $n,k,p,q$ and $(\text{diam}\Omega)^{1-n/p}||b/\rho^*_k(\mathcal A)||_{L^p(\Omega)}$.
\end{Theorem}

The case $k=n$ is the well known maximum principle, due originally to Aleksandrov for operators in the general form \eqref{1.1}; (see for example \cite{Aleksandrov 1960,Aleksandrov 1963,GTbook}).  

By adapting the proof of Corollary 4.2 in \cite{KT 2007}, we can obtain, from the proof of Theorem 1.1, the following variant, where the dependence on the diameter of $\Omega$ is replaced by that on the $L_q$ norm of  $\mathcal A / \rho^*_k(\mathcal A)$, thereby providing an extension to unbounded domains. To  fit better with later applications, we will express our assumed bounds on $\mathcal A$ in terms of $\mathcal T = \mathcal T(\mathcal A) := \text{trace} \mathcal A$.

\begin{Theorem} Assume additionally in Theorem 1.1, that  $k>n/2$ and $\mathcal A / \rho^*_k(\mathcal A) \in L^k(\Omega)$.
Then we have the estimate

\begin{equation} \label{1.9}
\sup_\Omega u \le \sup_{\partial\Omega} u + C R^{2-n/q} ||\frac{(Lu)^-}{\rho^*_k(\mathcal A)}||_{L^q(\Omega)},
\end{equation}

\noindent where $R = || \mathcal T/ \rho^*_k(\mathcal A)||^{k/n}_{L^k(\Omega)}$ and $C$ is a constant depending on $n,k,p$ and $R^{1-n/p}||b/\rho^*_k(\mathcal A)||_{L^p(\Omega)}$.  

\end{Theorem}

We will consider the cases, $k < n/2$ in Theorem 1.1 and $k\le n/2$ in Theorem 2.2, in conjunction with our treatment in Sections 2 and 3. We just remark here that Theorems 1.1 and 1.2 extend to these cases with $q >n/2$ and the conditions on the coefficient $b$ replaced by $b/\sqrt{\rho^*_k(\mathcal A) \lambda_{min}(\mathcal A)} \in L^p(\Omega)$ for $p>n$.  

As further consequences of Theorem 1.1, we also obtain, in Section 4, corresponding improvements of the local estimates for uniformly elliptic operators in \cite{KT 2007} which then provide further extensions of the fundamental  H\"older and Harnack estimates of Krylov and Safonov to the ellipticity cones  $\Gamma^*_k$. The reader is also referred to \cite{Krylov 2021}, \cite{Krylov 2021-2} and \cite{Safonov 2010} for extensions of the case $k=n$ to exponents $q\le n$, with $q$ at least close to $n$ and $p=n$. 

Finally, we remark here that Theorems 1.1 and 1.2, extend immediately to more general operators of the form,
\begin{equation}\label{1.10}
Lu =  a^{ij}D_{ij}u + b_iD_iu +cu,
\end{equation}
with scalar coefficient $c\le 0$ satisfying  $c/\rho^*(\mathcal A) \in L^q(\Omega)$.

\vspace{20mm}

\section{Proof of Theorem 1.1}\label{Section 2}

Theorem 1.1 is a consequence of a sharper version of the special case of operators $L_0$ of  the form, 
\begin{equation} \label{2.1}
L_0u:= \mathcal A. D^2 u,
\end{equation}
proved in Section 2 of \cite{KT 2007}, and the gradient estimate from Theorem 4.1 in \cite{TW 1999}. First we  recall the upper $k$-contact set of 
a function $u\in L^\infty(\Omega)$ is  defined by

\begin{equation}\label {2.2}
\Omega^+_k = \{x_0 \in \Omega |\  \exists \  k\text{-convex}\ v\in C^2(\Omega),\ \text{satisfying} \  u \le -v \ \text{in}\ \Omega, \ u(x_0) = - v(x_0)\}
\end{equation}

\noindent where a function $v \in C^2(\Omega)$ is called $k$-convex if $D^2 v \in \bar\Gamma_k(\Omega)$. Then from Section 2 of \cite{KT 2007}, it follows that if $u \in C^2(\Omega)\cap C^0(\bar\Omega)$ satisfies $u\le 0$ on $\partial\Omega$,  and  $\Omega_0$ is a domain containing 
$\Omega$, with $u$ extended to vanish on $\Omega_0 - \Omega$, so that $\Omega^+_{0,k} \subset \Omega$, we have the estimate,

\begin{equation} \label{2.3}
\sup_\Omega u \le C (\text{diam}\Omega_0)^{2-n/q}  ||\frac{(L_0u)^-}{\rho^*_k(\mathcal A)}||_{L^q(\Omega^+_{0,k})},
\end{equation}

\noindent where $C$ is a constant depending on $n,k$ and $q$. 

To proceed further we need the extended notion of  $k$-convex function introduced in \cite{TW 1999}, namely an upper semi-continous function $v$ in a domain $\Omega$ is called $k$-convex if any quadratic polynomial  $\varphi$ for which the difference $v-\varphi$ has a local maximum in $\Omega$ is $k$-convex. As in \cite{TW 1999} we denote the linear space of $k$-convex functions in $\Omega$ by $\Phi^k(\Omega)$. Equivalent characterisations are also provided in Lemmas 2.1 and 2.2 in \cite{TW 1999}. From Theorem 4.1 and Lemma 2.3 in \cite{TW 1999},  it then follows that $\Phi^k(\Omega)$ is contained in the Sobolev space $W^{1,r}_\text{loc}(\Omega)$ for $r < \frac{nk}{n-k}$ and moreover, we have the estimate: 

\begin{equation} \label{2.4}
|| Dv||_{L^r(\Omega^\prime)} \le C (\text{diam} \Omega) ^{(n-r)/r} \sup_\Omega |v|
\end{equation}

\noindent for any $v\in \Phi^k(\Omega)\cap L^{\infty}(\Omega)$ and subdomain $\Omega^\prime$ satisfying dist$(\Omega^\prime, \partial \Omega) \ge \kappa (\text{diam} \Omega)$, for some $\kappa \in (0,1/2)$, with constant $C$ depending on $n, k, r$ and $\kappa$. 

Now, defining the upper $k$-envelope of $u$ on $\Omega_0$ by

\begin{equation} \label{2.5}
u_0 = \inf \{ -v | v\in \Phi^k(\Omega_0), u\le -v \ \text{in}\  \Omega_0 \},
\end{equation}

\noindent we then have $u_0 \in \Phi^k(\Omega_0)$, $0 \le u \le u_0$ in $\Omega_0$ and $u_0 = u$ in $\Omega^+_{0,k}$. To apply the estimate
\eqref{2.3}, we suppose $\Omega$ lies in a ball of radius $R$  and take $\Omega_0$ to be the concentric ball of radius $2R$. Using the coordinate transformation $x \rightarrow (x-x_0)/R$, where $x_0$ denotes the centre of $\Omega_0$, we can then assume $R =1$. For convenience, by dividing $L$ by $\rho^*_k(\mathcal A)$, we can also assume $\rho^*_k(\mathcal A) = 1$. From the estimate \eqref{2.4}, applied to the function $u_0$ in $\Omega_0$, with  $\kappa = 1/4$ and $r$ chosen so that 
$$ \frac{1}{r} =  \frac{1}{q} - \frac{1}{p}, $$
we then have the estimate

\begin{equation} \label{2.6}
||b.Du||_{L^q(\Omega^+_{0,k})} \le C ||b||_{L^p(\Omega^+_{0,k})} \sup u, \\
\end{equation}

 \noindent where $C$ is a constant depending on $n,k,p$ and $q$. By combining \eqref{2.3} and \eqref{2.6}, we can obtain the desired estimate \eqref{1.6}, in the case when $||b||_p$ is sufficiently small. To proceed to the general case, we adapt the argument used by Pucci \cite{Pucci 1966 -1} in his treatment of the case $k=n$. Writing $M = \sup_\Omega u$, we let $\mathcal U_m = \mathcal U_{m,k}$, for $m\in [0,M)$ denote the $k$-upper contact set of the function $u-m$ in $\Omega_0$, so that $\mathcal U_0 = \Omega^+_{0,k}$ in our previous notation. For $m=M$, we set $\mathcal U_M = \{x\in \Omega \ | \ u=M\}$. For a sufficiently large positive integer $N$, we then let $m_0, \dots, m_N$ be an increasing sequence in $[0,M]$, such that $m_0 = 0$, $m_N = M$ and
 
 \begin{equation}\label{2.7}
 \int _{\mathcal U_{i-1} - \mathcal U_i} |b|^p = \frac{1}{N^{p/q}}  \int _{\mathcal U_0 - \mathcal U_M} |b|^p
 \end{equation}
 
 \noindent for $i=1, \dots,N$, and $\mathcal U_i = \mathcal U_{m_i}$. Now we define
 
$$ y_i := \int_{\mathcal U_{i-1} - \mathcal U_{i}} |b.Du|^q,$$

 \noindent so that
  
 $$ \Sigma^{N}_{j=i} y_j= \int_{\mathcal U_{i-1}} |b.Du|^q $$
 
\noindent  for $i = 1, \dots, N$, since $Du=0$ a.e on $\mathcal U_N$.  Applying the estimates, \eqref{2.3} and \eqref{2.6}, to the functions $u-m_i$ and using \eqref {2.7}, we then obtain, for $i=1, \dots,N$, the estimates
 
 \begin{equation}\label{2.8}
 y_i \le \frac{C_0}{N}  ||b||^q_{L^p(\mathcal U_0)} \{ \Sigma^N_{j=i} y_j + ||(Lu)^-||^q_{L^q(\mathcal U_0)}\},
 \end{equation}
 
 \noindent for a further constant $C_0$ depending on $n, k, p$ and $q$.  Now setting 
 
 $$\theta = C_0 ||b||^q_{L^p(\mathcal U_0)}, \ \alpha = \theta/ (N-\theta), \  \mu = ||(Lu)^-||^q_{L^q(\mathcal U_0)}$$
 
\noindent and taking $N> \theta$, we can rewrite \eqref{2.8} as 

$$ y_i \le \alpha( \Sigma^N_{j=i+1} y_j + \mu)$$

\noindent In accordance with the discrete Gronwall inequality, we then have, by iteration, 
$$y_i \le \alpha (1+\alpha)^{(N-i)} \mu,$$

\noindent for $i=1, \dots,N$. Summing over $i$, we now obtain 

$$ \int_{\mathcal U_0} |b.Du|^q \le \{(1+\alpha)^N-1\} \mu = \{(\frac{N}{N-\theta})^N -1 \}\mu $$


\noindent so that letting $N\rightarrow \infty$, we obtain the estimate

\begin{equation}\label{2.9}
 ||b.Du||_{L^q(\mathcal U_0)} \le (e^{\theta} -1)^{1/q} ||(Lu)^-||_{L^q(\mathcal U_0)}.
\end{equation}

\noindent Note that we can also infer \eqref{2.9} directly from \eqref{2.8} and the Lemma in Section 3 of \cite {Pucci 1966 -1}, which corresponds to a cruder version of our preceding estimate with finite $N>\theta$.

Combining \eqref{2.3} and \eqref{2.9} and returning to our original coordinates, with $u$ replaced by $u-\sup_{\partial\Omega}u$, we obtain the following sharper version of our desired estimate \eqref{1.7}  for functions $u \in C^2(\Omega)\cap C^0(\bar\Omega)$, 
\begin{equation}\label{2.10}
\sup_\Omega u \le \sup_{\partial\Omega}u \ + \ Cd^{2-n/q} \exp(C_0 d^{(1-n/p)q}||b/\rho^*_k(\mathcal A))||^q_{L^p(\Omega^+_k)}) ||\frac{(Lu)^-}{\rho^*_k(\mathcal A)}||_{L^q(\Omega^+_k)},
\end{equation}

\noindent where $d = \text{diam} \Omega$, $C$ and $C_0$ are constants depending respectively on $n,k,q$ and $n,k,p,q$. By following the approximation argument given, for example, in the proof of the case $k=n$ in \cite{GTbook}, we may extend the estimate to functions $u \in W^{2,q}_{loc} (\Omega)\cap C^0(\bar\Omega)$, thereby completing the proof of Theorem 1.1.

We can extend the estimate \eqref{2.10} to the cases $k<n/2$ through a completely different approach by adapting that used for divergence structure operators in \cite{T 1971}. Under our previous normalisations, $R= \rho^*_k(\mathcal A) = 1$, we consider the function $w$ given by

\begin{equation}\label{2.11}
w= \log \frac{M+\mu}{M-u+\mu},
\end{equation}

\noindent where $\mu = ||(Lu)^-||_{L^q(\Omega)}$ for some $q>n/2$, so that

\begin{align*}
 L_0 w &= \mathcal A Dw.Dw - b.Dw + \frac{1}{M-u+\mu} Lu \\
        & \ge - \frac{1}{M-u+\mu} (Lu)^-  -\frac{|b|^2}{4\lambda_{min}(\mathcal A)}.
 \end{align*}    
 
 \noindent Applying the estimate \eqref{2.3}, with $\Omega_0 = \Omega$, we then obtain
 
$$ \sup w \le C(1+ ||b^2 /4\lambda_{min}(\mathcal A)||_{L^q(\Omega)})$$ 

\noindent and hence, from  \eqref{2.11}  and taking account of our normalisations, we obtain, in place of \eqref{2.10},

\begin{equation}\label{2.12}
\sup_\Omega u \le \sup_{\partial\Omega}u \ + \ d^{2-n/q} \exp C\{1+d^{2q-n)}||b^2/\rho^*_k(\mathcal A)\lambda_{min}(\mathcal A)||^q_{L^q(\Omega}\} ||\frac{(Lu)^-}{\rho^*_k(\mathcal A)}||_{L^q(\Omega)},
\end{equation}

\noindent for all $k$ and any $q>n/2$, where $C$ is a constant depending on $n$, $k$, and $q$.

\vspace{20mm}

\section{Proof of Theorem 1.2 }\label{Section 3}

Unlike the case $L = L_0$ in \cite{KT 2007}, Theorem 1.2 does not immediately follow from the corresponding local maximum principle. However the approach is still similar and we can combine the proofs, with the resultant local maximum principle formulated in the next section as Theorem 3.1.
Accordingly we fix a ball $B: = B_R: = B_R(y)$  in $\mathbb R^n$ of centre $y$ and radius $R$, such that the intersection  $\Omega\cap B$ is non-empty. Under the hypotheses of Theorem 1.2 and our previous normalisations, $B = B_1(0)$ and $\rho^*_k(\mathcal A) = 1$, we define a cut-off function $\eta$ by

\begin{equation}\label{3.1}
\eta = [(1 - |x|^2 )^+]^{\beta}
\end{equation}

\noindent for some $\beta \ge 1$, to be chosen. Then setting $v = \eta (u^{+})^2$, for $u \in C^2(\Omega)\cap C^0(\bar\Omega)$ satisfying $u\le 0$ on $\partial\Omega\cap B$, we now compute, in $\Omega\cap B\cap \{u>0\}$,

\begin{equation}\label{3.2}
\begin{aligned} 
L_0 v: = a^{ij}D_{ij} v &=  u^2a^{ij} D_{ij} \eta  + 4u a^{ij} D_i \eta
D_ju + 2\eta a^{ij} D_i uD_ju + 2\eta u L_0u \\
&\geq - C\mathcal T \eta^{1 - 2/\beta}u^2 - \eta b.Du^2  + 2\eta u Lu,
\end{aligned}
\end{equation}
\noindent  where $C$ is a constant depending on $n$ and $\beta$. Now we fix the point $x_0$ so that $u(x_0) = \sup_\Omega u = \sup_\Omega v$ and take $\beta=2$. Applying the estimate \eqref{2.3} with $\Omega_0 = B_2$, we then obtain

\begin{equation}\label{3.3}
\sup_\Omega v \le C || (\eta^{1 - 2/\beta} \mathcal Tu^2+ \eta (b.Du^2) ^- + \eta u (Lu)^-||_{L^q(\mathcal U_0)}
\end{equation}

\noindent for a constant $C$ depending on $n,k,q$ and $\beta$, where now $\mathcal U_0 \subset B \cap \Omega$ denotes the upper $k$-contact set of the function $v$ in $\Omega_0$.

\noindent To prove the global estimate, \eqref {1.8}, we now fix the point $x_0$ so that $u(x_0) = \sup_\Omega u = \sup_\Omega v$ and take $\beta=2$ in \eqref{3.3}.  Consequently, if 

\begin{equation} \label{3.4}
||\mathcal T||_{L^q(\Omega)} \le 1/2C,
\end{equation}
where $C$ is the constant in \eqref{3.3}, with $\beta=2$, we have 

$$\sup_\Omega v \le 2C|| \eta (b.Du^2) ^- + \eta u (Lu)^-||_{L^q(\mathcal U_0)}$$

\noindent so that we can adapt the proof of Theorem 1.1 to handle the coefficient $b$. Here we let $\mathcal U_m = \mathcal U_{m,k}$, for $m\in [0,M)$, denote the $k$-upper contact set of the function $v_m:= \eta(u^2-m^2)$ in $\Omega_0$ and define 

$$ y_i := \int_{\mathcal U_{i-1} - \mathcal U_{i}}(\eta |b.Du^2|)^q,$$

 \noindent so that now

$$ \Sigma^{N}_{j=i} y_j = \int_{\mathcal U_{i-1}} (\eta |b.Du^2|)^q $$
 \noindent  for $i = 1, \dots, N$. 
 
 To proceed further we also need to estimate $D\eta/\eta$ in $\mathcal U_m$. Letting $v_0$ denote the upper $k$-envelope of $v_m$ in $\Omega_0$, it follows that 
 
$$ (1-(|x|/2)^{2-n/k} \sup v_m  \le  v_0 \le \eta \le (1-|x|^2) \sup v_m$$                                                                                             
 
 \noindent in $\mathcal U_m$, which implies 
 
 $$ |D\eta|/\eta \le \frac{\beta}{ 1-(1/2)^{(2k/n) -1}}:= \gamma $$

\noindent  in $\mathcal U_m$. 
 
 \noindent Now denoting $v^i = v_{m_{i-1}}$, we then have, in the upper contact set $\mathcal U_{i-1}$, 
 
 \begin{align*}
  \eta |b.Du^2| & \le |b| \{ |Dv^i| + |D\eta|(u^2 - m^2_{i-1})\} \\
                      & \le |b| \{ |Dv^i| +\gamma \sup v^i\},
 \end{align*}
so that using the estimates \eqref{2.3} and \eqref{2.6} we obtain again the fundamental inequality \eqref {2.8}, for a constant $C$ depending on $n,k,p, q$ and $\beta$. We remark that the case $p=k=n=q$ is much simpler here in that we can just use concavity of the $k$-upper envelope of $v^i$ to estimate $|Dv^i| \le v^i$ in $\mathcal U_{i-1}$. The proof can then be further simplified by taking $v=\eta u^+$. 

To conclude the estimate \eqref{1.8}, we then need to fix $\beta=2$ and an appropriate radius $R$ satisfying the condition \eqref{3.4} in our original coordinates which is done by taking 

\begin{equation} \label {3.5}
R = C|| \mathcal T/\rho^*_k(\mathcal A)||^{q/n}_{L^q(\Omega)} 
\end{equation}
for a further constant $C$ depending on $n,k$ and $q$. 

Accordingly we obtain an estimate of the form \eqref{2.10}, with the diameter $d$ replaced by $|| \mathcal T/\rho^*_k(\mathcal A)||^{q/n}_{L^q(\Omega)}$, thereby completing the proof of Theorem 1.2.

In order to embrace the cases $k<n/2$, we revisit the corresponding extension  of Theorem 1.1 and set $v=\eta w$, where $w$ is given by \eqref{2.11}
and as above, the cut-off function $\eta$ is given by \eqref{3.1}, with $\beta = 2$. Then we obtain, in place of \eqref{3.2},

\begin{equation}\label{3.6}
\begin{aligned} 
L_0 v &=  wa^{ij} D_{ij} \eta  + 2 a^{ij} D_i \eta D_jw + \eta a^{ij} D_i wD_jw -\eta b.Dw + \frac{\eta}{M-u+\mu} Lu \\
          &\geq - C\mathcal T w  -\frac{|b|^2\eta}{4\lambda_{min}(\mathcal A)}  - \frac{\eta}{M-u+\mu} (Lu)^- .
\end{aligned}
\end{equation}

Now choosing $R$ as in \eqref{3.5} and applying the estimate \eqref{2.3} in $B_{2R}$, we obtain an estimate of the form \eqref{2.12} with $d$ replaced by $|| \mathcal T/\rho^*_k(\mathcal A)||^{q/n}_{L^q(\Omega)}$. 


\vspace{20mm}

\section{Local estimates }\label{Section 4}

We begin with an extension of the local maximum principle in Theorem 4.1 in \cite{KT 2007}, under a slight sharpening of the hypotheses of Theorem 1.2.

\begin{Theorem} Assume additionally in Theorem 1.1, that  $\mathcal A / \rho^*_k(\mathcal A) \in L^r(\Omega)$ for $r > \max\{k,n/2\}$. Then for any ball
$B = B_R = B_R(y)$ intersecting $\Omega$, concentric ball $B_{\sigma R}$, $0<\sigma<1$, and $\kappa >0$ we have the estimate

\begin{equation} \label{4.1}
\sup_{\Omega\cap B_\sigma} u \le \sup_{\partial\Omega \cap B} u + C \big\{ [R^{-n} \int_{\Omega \cap B} (u^+)^\kappa \ ]^{1/\kappa}  + R^{2-n/q}    ||\frac{(Lu)^-}{\rho^*_k(\mathcal A)}||_{L^q(\Omega\cap B)} \big\},
\end{equation}

\noindent where $C$ is a constant depending on $n,k,p,q,r, \sigma, \kappa$, $R^{1-n/p}||b/\rho^*_k(\mathcal A)||_{L^p(\Omega\cap B)}$ and $R^{-n/r}||\mathcal T/\rho^*_k(\mathcal A)||_{L^r(\Omega \cap B)}$ 

\noindent If $k < n/2$, the estimate \eqref{4.1} holds with $b/\rho^*_k$ replaced by $b/\sqrt{\rho^*_k \lambda_{min}}$.

\end{Theorem}

To obtain the estimate \eqref{4.1}, in the cases $k\ge n/2$, we can simply return to the proof of Theorem 1.2, in Section 3, and not use condition \eqref {3.4} in the estimate \eqref{3.3}. We then obtain the estimate

\begin{equation} \label{4.2}
\sup_ {\Omega\cap B} v \le C \exp(C_0 ||b||^q_{L^p(\Omega\cap B)}) || \eta^{1 - 2/\beta}\mathcal T(u^+)^2 +\eta u^{+} (Lu)^-||_{L^q(\Omega\cap B)}, 
\end{equation}

\noindent where, corresponding  to \eqref{2.10}, $C$ and $C_0$ are constants depending respectively on $n,k,\beta$ and $n,k,\beta, p$. Since

$$|| \eta^{1 - 2/\beta}\mathcal T(u^+)^2 +\eta u^{+} (Lu)^-||_q \le ||\mathcal T||_r ||v^{1-2/\beta}(u^+)^{4/\beta}||_{r^*} + ||v^{1/2}(Lu)^- ||_q, $$

\noindent where 

$$ \frac{1}{r^*} = \frac{1}{q} - \frac{1}{r} >0, $$

\noindent we can then infer \eqref {4.1}, by taking $\beta= 4r^* /\kappa$. The cases $k<n/2$ are essentially already proved in \cite{KT 2007}, as we can estimate from \eqref{3.2},

$$L_0 v \geq - C(\mathcal T +\frac{|b|^2}{\lambda_{min}(\mathcal A)} ) \eta^{1 - 2/\beta}u^2 + 2\eta u Lu $$

\noindent where $C$ depends on $n$ and $\beta$. Instead of \eqref{4.2}, we now obtain the estimate

$$\sup_ {\Omega\cap B} v \le C || \eta^{1 - 2/\beta}(\mathcal T + \frac{|b|^2}{\lambda_{min}(\mathcal A)})(u^+)^2 +\eta u^{+} (Lu)^-||_{L^q(\Omega\cap B)}, $$

\noindent where $C$ depends on $n, k, q,\beta$. Fixing $p=2r >2q$ and $r^*$ as before, we then infer the cases  $k < n/2$ of Theorem 4.1.

From Theorem 1.1 also follow corresponding extensions of the H\"older and Harnack estimates in Theorem 4.3 in \cite{KT 2007} for uniformly elliptic operators, where the condition on $\mathcal A$ in Theorem 3.1 is strengthened to  $\mathcal A / \rho^*_k(\mathcal A) \in L^\infty(\Omega)$.
Defining 

$$ a_k := \sup_\Omega \frac{\lambda_{max}} {\rho^*_k}(\mathcal A), $$

\noindent we then have from inequality \eqref{1.6}, the uniformly elliptic condition,

$$\frac{\lambda_{max}}{\lambda_{min}}(\mathcal A) \le (\frac{n}{k})^k {\binom{n}{k}}^{-1}  a^k_k \le a^k_k,$$

\noindent which is a refinement of inequality (4.3) in \cite{KT 2007}.

Our H\"older and Harnack estimates now follow as corollaries of the following weak Harnack inequality for supersolutions, which was proved in 
\cite{Trudinger 1980} in the case $k=n, q=2n$, by adapting the key ideas in \cite{KS 1980}. 

\begin{Theorem} Assume additionally to the conditions on $L$ in Theorem 1.1, that  $L$ is uniformly elliptic in $\Omega$. Then, for any ball
$B = B_R = B_R(y)\subset \Omega$, $u \ge 0, \in W^{2,q}(B)$, for $q=k$ if $k >n/2$, $q>n/2$ if $k\le n/2$, and $0<\sigma, \tau <1$ , there exists $\kappa >0$ such that
 
\begin{equation} \label{4.3} 
 [R^{-n} \int_{B_{\sigma R}} u^\kappa \ ]^{1/\kappa}  \le C \big\{ \inf_{B_{\tau R}}u  + R^{2-n/q}  ||\frac{(Lu)^+}{\rho^*_k(\mathcal A)}||_{L^q(B_R)} \big\}
\end{equation}
where $\kappa$ and $C$ depend on $n, k, p, q, \sigma, \tau, a_k$ and $R^{1-n/p}||b/\rho^*_k(\mathcal A)||_{L^p(\Omega}$.
\end{Theorem}

\begin{Corollary}
Assume that L satisfies the hypotheses of Theorem 4.2 and $u\in  W^{2,q}(B)$. Then for any concentric ball $B_\sigma = B_{\sigma R}(y)$,  $0<\sigma <1$, we have the oscillation estimate

\begin{equation}\label{4.4}
\underset{B_\sigma}{\text{osc}}\ u \le C\sigma^{\alpha} \big\{ \underset{B}{\text{osc}}\ u+ R^{2-n/q}  ||\frac{Lu}{\rho^*_k(\mathcal A)}||_{L^q(B_R)} \big\}
\end{equation}
where $\alpha>0$ and $C$ depend on  $n, k, p, q,  a_k$ and $R^{1-n/p}||b\rho^*_k(\mathcal A)||_{L^p( B)}$. Furthermore if $u\ge 0$, we have the Harnack inequality

\begin{equation}\label{4.5}
\sup_{B_\sigma} u \le C \big\{ \inf_ {B_\sigma} u  + R^{2-n/q}  ||\frac{Lu}{\rho^*_k(\mathcal A)}||_{L^q(B_R)} \big\}
\end{equation}
where $C$ depends on the same quantities as in \eqref{4.3}.

\end{Corollary}

Note that the case $p=k=q=n$ is proved in \cite{Safonov 2010} and its derivation from the corresponding case in Theorem 1.1 is  more technically delicate than the cases where $p>n$. The cases $k\le n/2$, where we can take $p=2q$,  are essentially already covered in \cite{KT 2007}. Here we will just indicate how the cases $p>n$ follow by simple modifications of previous approaches, when $k=n$. For convenience we will just adapt the presentation in \cite{Tokyo lectures}, although, as remarked above the main ideas go back to \cite{KS 1980}. The first step in the proof of Theorem 4.2 corresponds to a weak form of the estimate \eqref{4.3}, where $\sigma < \tau$ and the left hand side is replaced by $\inf_{B_{\sigma R}}u$. Normalising 
$R= \lambda_{min}(\mathcal A) = 1$, we define a comparison function $v$ by

$$v(x) := \frac{|x|^{-\beta} - 1}{\sigma^{-\beta} - 1} \inf_{B_\sigma} u,$$

\noindent for a constant $\beta \ge \sup \mathcal T -2$. We then obtain

$$ Lv \ge -\frac{\beta|x|^{-\beta-2}}{\sigma^{-\beta}-1}(b.x)\inf_{B_\sigma} u,$$

\noindent in $B-B_\sigma$, so that by applying Theorem 1.1 to the function $v-u$, and assuming $\sigma \le 1/2$, we obtain

$$ u(x)-v(x) \ge - C || \frac {|b|}{\sigma} \inf_{B_\sigma} u+ (Lu)^+|| _{L^q(B-B_{\sigma})}, $$

\noindent for $\sigma<|x|<1$, where $C$ depends on $n, k, p, q$, $\sup \mathcal T$ and $||b||_{L^p( B)}$.

By removing the normalisation, $R=1$,  and taking account of the more explicit estimate \eqref{2.10}, we then obtain an estimate

\begin{equation}\label{4.6}
 \inf_{B_{\sigma R}}u \le C \big\{ \inf_{B_{\tau R}}u  + R^{2-n/p}  ||(Lu)^+||_{L^p(B_R)} \big\},
\end{equation}

\noindent provided $R$ is sufficiently small, in the sense that

$$R^{1-n/p}||b||_{L^p( B_R)} \le 1/C,$$

\noindent where now $C$ depends on  $n, k, p, q, \sigma,\tau$ and $a_k$. 

The extension to general $R$ can then be achieved by standard covering arguments. In particular, after normalising $R=1$ as before, we can apply the estimate \eqref{4.6}, with $\sigma =1/3, \tau = 2/3$, to a finite chain of overlapping balls $B_\rho(y_i)$ for $i=1,\ldots N$ satisfying 

$$ \rho^{1-n/p}||b||_{L^p( B)} \le 1/C,$$

\noindent whose centres $y_i$  are equally spaced along a line segment  joining  $y_1 \in B_{\sigma - \rho/3}$ to $y_N \in B_\tau $ satisfying
$|y_{i+1} -y_i| =\rho/3$, $i=1, \dots N-1$, $|y_N - y_0| \in (\rho/3, 2\rho/3)$, where $y_0 \in B_\tau$ satisfies $u(y_0) = \inf_{B_\tau}u$. By iteration, we then obtain an estimate \eqref {4.6}, under the hypotheses of Theorem 4.2, (with $\lambda_{min}(\mathcal A) = 1$), where now the constant $C$ depends additionally on $R^{1-n/p}||b||_{L^p(B_R)}$.

The weak Harnack inequality, Theorem 4.2  then follows by combining the estimate (4.6) with the local maximum principle in Theorem 4.1 applied to the function $1-u$,  and using the key measure theoretic argument of Krylov and Safonov, as presented for example in \cite{Tokyo lectures}.  

The results of this section also extend to operators of the form \eqref{1.10} with $c/\rho^*(\mathcal A) \in L^q(\Omega)$

\vspace{20mm}

\section{Application to the Pucci conjecture }\label{Section 5}

We conclude this paper by indicating the application to the corresponding improvement of the condition on the drift term in \cite{T 2020} in the case $q=k > n/2$.  First we note that, by defining a uniform ellipticity constant  $a_0$ for the operator $L$ by

$$ a_0 = \sup_{\Omega} \frac{\mathcal T}{\lambda_{min}}(\mathcal A), $$

\noindent it follows from  Section 3 in \cite{T 2020}, that $\mathcal A \in \Gamma^*_k$ if 

\begin{equation} \label{5.1}
 \chi := k -n[1- (n-1) /a_0] >0. 
\end{equation}

\noindent Note here that we are not normalising $\mathcal T = 1$ as in \cite{T 2020} and $a_0 = 1/\delta$ where $\delta$ is the ellipticity constant used in \cite{Pucci 1966} and \cite{T 2020}.

Moreover from the first inequalities on page 116 there, with $k$ replaced by $k-\chi$, we can obtain an explicit estimate from below for  $\rho^*_k( \mathcal A)$. As there, we first note that for  $\lambda_i > 0$, $i = 1, \ldots, n$, the inner product $\lambda .\mu$ does not increase when we rearrange the components of $\lambda$ in increasing order and those of $\mu$ in decreasing order. Consequently if $\lambda = \lambda(\mathcal A)$ and $\mu \in \Gamma_k$ are so ordered, with $\mu_n \le 0$, we obtain

\begin{align*}
\lambda.\mu &\ge \lambda_1\ \{ \sum_{i<n} \mu_i \ +  [a_0-(n-1)] \mu_n \}\\
                     & \ge \lambda_1 \ \{ \sum_{i<n} \mu_i \ + \frac{(k-\chi)(n-1)}{(n-k+\chi)}\ \mu_n \} \\
                     & = \lambda_1 \ \{S_1(\mu) \ + \frac{n(k-1-\chi)}{n-k+\chi} \ \mu_n \}.
\end{align*}

\noindent Now, using the fundamental inequality,

\begin{equation} \label{5.2}
(n-k)S_1(\mu) \ + n(k-1)\mu_n \  > 0
\end{equation}
for $\mu \in \Gamma_k$, we then obtain, for $k>1$,

\begin{align*}
\lambda.\mu & \ge \lambda_1S_1(\mu) \frac{(n-1)\chi}{(n-k+\chi)}\\
                     & = \lambda_1S_1(\mu)\frac{\chi a_0}{n(k-1)},
 \end{align*}  
 and hence our desired estimate,
 
 \begin{equation}\label{5.3}
 \rho^*_k( \mathcal A) \ge \frac{\chi a_0}{n(k-1)}\lambda_{min},
 \end{equation}
 for $k>1$. Note that we can cover the case when $\mu_n > 0$ by assuming $\chi a_0 \le (k-1)/(n-1)$ and the argument above extends to the case $k = 1$, with $\chi = 0$,  and is consistent with the equality $\rho^*_1 = \lambda_{min}$.

Consequently Theorems 1.1,1.2, together with their extensions to $1< k \le n/2$, and Theorem 4.1, for $k>1$, extend to uniformly elliptic operators satisfying \eqref{5.1}, in place of 
$\mathcal A \in \Gamma^*_k$, with $\rho_k^*(\mathcal A)$ replaced by $ {\lambda_{min}}(\mathcal A)\chi a_0$. Similarly, Theorem 41 and Corollary 4.2 extend with additionally $a_k$ replaced by $1/\chi$.



\vspace{20mm}

\baselineskip=12.8pt
\parskip=3pt
\renewcommand{\baselinestretch}{1.38}


\begin{thebibliography}{99999}

\bibitem{Aleksandrov 1960} A.D.Alexandrov, Certain estimates for the Dirichlet problem, Dokl. Akad. Nauk. SSSR, 134, (1960), 1001-1004,
( Russian). English translation : Soviet Math. Dokl.{\bf 1}(1960),1151-1154.

\bibitem{Aleksandrov 1963} A.D.Alexandrov, The method of normal map in uniqueness problems and estimations for elliptic equations, in: Seminari dell'Istituto Nazionale di Alta Matematica 1962-63, Ediz. Cremonese, Rome, (1965), 744-786.








\bibitem{GTbook} D. Gilbarg and N.S. Trudinger, Elliptic partial differential equations of the second order, Second edition, 
Springer-Verlag, Berlin, 1983 (reprinted 2001).

\bibitem{Krylov 2021}N.V. Krylov, On stochastic equations with drift in $L_d$. Ann. Probab. 49 (2021), 2371-2398. 


\bibitem{Krylov 2021-2}N.V. Krylov, A review of some new results in the theory of linear elliptic equations with drift in $L_d$. Anal. Math. Phys. 11 (2021),  Paper No. 73, 13 pp. 




\bibitem{KS 1980} N.V.Krylov and M. V. Safonov, A certain property of solutions of parabolic equations with
measurable coefficients, Izv. Akad. Nauk SSSR Ser. Mat. 44 (1980), 161-175, (Russian); English translation: Math. USSR Izv.16 (1981), 151-164.

\bibitem{KT 2007} H-J Kuo and N.S. Trudinger, New maximum principles for linear elliptic equations, Indiana Univ. Math. J. 56 (2007), 2439-2452.



\bibitem{Pucci 1966} C. Pucci, Operatori ellitiche estremanti, Ann. Mat. Pura App. 72, (1966), 141-170.

\bibitem{Pucci 1966 -1} C. Pucci, Limitazione per soluzione di equazione ellitiche, Ann. Mat. Pura App. 74, (1966), 15-30


\bibitem{Safonov 2010} M. V. Safonov,  Non-divergence Elliptic Equations of Second Order with Unbounded Drift, Nonlinear partial differential equations and related topics, Amer. Math. Soc. Transl. Ser. 2, 229, Adv. Math. Sci., 64, Amer. Math. Soc., Providence, RI, 2010, 211-232, 

\bibitem{T 1971} N.S.Trudinger, Linear elliptic operators with measurable coefficients,  Ann. Scuola Norm. Sup. Pisa Cl. Sci. (3) 27 (1973), 265-308.

\bibitem{Trudinger 1980} N.S.Trudinger, Local estimates for subsolutions and supersolutions of general second order elliptic quasilinear equations, Invent. Math., 61 (1980), 67-79.

\bibitem{Tokyo lectures} N.S. Trudinger, Lectures on Nonlinear Elliptic Equations of Second Order, Lectures in Mathematical Sciences, University of Tokyo, New Series Vol. 9, 1995.

\bibitem{Trudinger 1995} N. S. Trudinger, On the Dirichlet problem for Hessian equations, Acta Math. 175 (1995), 151-164.

\bibitem{T 2020} N. S. Trudinger, Remarks on the Pucci conjecture,  Indiana Univ. Math. J. 69 (2020), 109- 118.

\bibitem {TW 1997} N.S. Trudinger and X.-J Wang, Hessian measures. I, Topol. Methods Nonlinear Anal. 10 (1997), 225-239, Dedicated to Olga Ladyzhenskaya.

\bibitem {TW 1999} N.S. Trudinger and X.-J Wang, Hessian measures. II, Ann. of Math. (2) 150 (1999), 579-604.







\end{thebibliography}
\end{document}